\documentclass{amsart}
\newtheorem{theorem}{Theorem}[section]
\newtheorem{lemma}[theorem]{Lemma}

\theoremstyle{definition}
\newtheorem{definition}[theorem]{Definition}
\newtheorem{example}[theorem]{Example}
\theoremstyle{remark}
\newtheorem{remark}[theorem]{Remark}

\numberwithin{equation}{section}
\begin{document}
%%%%%%%%%%%%%%%%%%%%%%%%%%%%%%%%%%%%%%%%%%%%%%%%%%%%%%%%%%%%%%%%%%%%%%%%%%%%

\title{Generalizations of the Lax-Milgram theorem}
%%%%%%%%%%%%%%%%%%%%%%%%%%%%%%%%%%%%%%%%%%%%%%%%%%%%%%%%%%%%%%%%%%%%%%%%%%%%

\author{Dimosthenis Drivaliaris}
\address{Department of Financial and Management Engineering\\
University of the Aegean
31, Fostini Str.\\
82100 Chios\\
Greece}
\email{d.drivaliaris@fme.aegean.gr}
\author{Nikos Yannakakis}
\address{Department of Mathematics\\
School of Applied Mathematics and Natural Sciences\\
National Technical University of Athens\\
Iroon Polytexneiou 9\\
15780 Zografou\\
Greece}
\email{nyian@math.ntua.gr}
\commby{}
%%%%%%%%%%%%%%%%%%%%%%%%%%%%%%%%%%%%%%%%%%%%%%%%%%%%%%%%%%%%%%%%%%%%%%%%%%%%

\begin{abstract}
We prove a linear and a nonlinear generalization of the Lax-Milgram theorem. In particular we give sufficient conditions for a real-valued function defined on the product of a reflexive Banach space and a normed space to represent all bounded linear functionals of the latter. We also give two applications to singular differential equations. 
\end{abstract}

\maketitle
%%%%%%%%%%%%%%%%%%%%%%%%%%%%%%%%%%%%%%%%%%%%%%%%%%%%%%%%%%%%%%%%%%%%%%%%

\section{Introduction}
%%%%%%%%%%%%%%%%%%%%%%%%%%%%%%%%%%%%%%%%%%%%%%%%%%%%%%%%%%%%%%%%%%%%%%%%
%%%%%%%%%%%%%%%%%%%%%%%%%%%%%%%%%%%%%%%%%%%%%%%%%%%%%%%%%%%%%%%%%%%%%%%%

The following generalization of the Lax-Milgram Theorem was proved recently by
An, Du, Duc and Tuoc in \cite{An}. 
%%%%%%%%%%%%%%%%%%%%%%%%%%%%%%%%%%%%%%%%%%%%%%%%%%%%%%%%%%%%%%%%%%%%%%%%

\begin{theorem}
\label{theorem1}
Let $X$ be a reflexive Banach space over $\mathbb R$, $\{ X_{n}\}_{n\in\mathbb N}$ be an increasing sequence of closed subspaces of $X$ and $V=\bigcup _{n\in\mathbb N}X_{n}$. Suppose that
\begin{equation*}
A:X\times V\longrightarrow \mathbb R
\end{equation*}
is a real-valued function on $X\times V$ for which the following hold:
\begin{itemize}
\item[\textit{(a)}] $A_{n}=A\left| _{X_{n}\times X_{n}}\right.$ is a bounded bilinear form, for all $n\in\mathbb N$.
\item[\textit{(b)}] $A(\cdot ,v)$ is a bounded linear functional on $X$, for all $v\in V$.
\item[\textit{(c)}] $A$ is coercive on $V$, i.e. there exists $c>0$ such that
\begin{equation*}
A(v,v)\geq c\| v\|^{2},
\end{equation*}
for all $v\in V$.
\end{itemize}
Then, for each bounded linear functional $v^*$ on $V$, there exists $x\in X$ such that 
\begin{equation}
\nonumber
A(x,v)=\langle v^*,v\rangle ,
\end{equation}
for all $v\in V$.
\end{theorem}
%%%%%%%%%%%%%%%%%%%%%%%%%%%%%%%%%%%%%%%%%%%%%%%%%%%%%%%%%%%%%%%%%%%%%%%%

In this paper our aim is to prove a linear and a nonlinear extension of Theorem \ref{theorem1}. In the linear case we use a variant of a theorem due to Hayden \cite{Hayden1, Hayden2} and thus manage to substitute the coercivity condition in (c) of the previous theorem with a more general $\text{inf-sup}$ condition. In the nonlinear case we appropriately modify the notion of type $M$ operator and use a surjectivity result for monotone, hemicontinuous, coercive operators. We also present two examples to illustrate the applicability of our results.
%%%%%%%%%%%%%%%%%%%%%%%%%%%%%%%%%%%%%%%%%%%%%%%%%%%%%%%%%%%%%%%%%%%%%%%%

All Banach spaces considered are over $\mathbb R$. Given a Banach space $X$, $X^{*}$ will denote its dual and $\langle\cdot ,\cdot\rangle$ their duality product. Moreover if $M$ is a subset of $X$, then $M^\perp$ will denote its annihilator in $X^\ast$ and if $N$ is a subset of $X^\ast$, then $^\perp N$ will denote its preannihilator in $X$.
%%%%%%%%%%%%%%%%%%%%%%%%%%%%%%%%%%%%%%%%%%%%%%%%%%%%%%%%%%%%%%%%%%%%%%%%

\section{The linear case}
%%%%%%%%%%%%%%%%%%%%%%%%%%%%%%%%%%%%%%%%%%%%%%%%%%%%%%%%%%%%%%%%%%%%%%%%

To prove our main result for the linear case we need the following lemma which is a variant of \cite[Theorem 12]{Hayden1} and \cite[Theorem 1]{Hayden2}:
%%%%%%%%%%%%%%%%%%%%%%%%%%%%%%%%%%%%%%%%%%%%%%%%%%%%%%%%%%%%%%%%%%%%%%%%

\begin{lemma}
\label{lemma1}
Let $X$ be a reflexive Banach space, $Y$ be a Banach space and
\begin{equation*}
A:X\times Y\longrightarrow \mathbb R
\end{equation*}
be a bounded, bilinear form satisfying the following two conditions:
\begin{itemize}
\item[\textit{(a)}] $A$ is non-degenerate with respect to the second variable, i.e., for each\\
$y\in Y\setminus\{ 0\}$, there exists $x\in X$ with $A(x,y)\ne 0$.
\item[\textit{(b)}] There exists $c>0$ such that
\begin{equation*}
\sup_{\| y\| =1}|A(x,y)|\geq c\| x\| ,
\end{equation*}
for all $x\in X$.
\end{itemize}
Then, for every $y^{*}\in Y^{*}$, there exists a unique $x\in X$ with
\begin{equation*}
A(x,y)=\langle y^{*},y\rangle ,
\end{equation*}
for all $y\in Y$.
\end{lemma}
%%%%%%%%%%%%%%%%%%%%%%%%%%%%%%%%%%%%%%%%%%%%%%%%%%%%%%%%%%%%%%%%%%%%%%%%
\begin{proof}
Let $T:X\longrightarrow Y^{*}$ with $\langle Tx,y\rangle =A(x,y)$,
for all $x\in X$ and all $y\in Y$. Obviously $T$ is a bounded, linear map. Since, by (b), $\| Tx\|\geq c\| x\|$, for all $x\in X$, $T$ is one-to-one. To complete the proof we need to show that $T$ is onto. 

Since $A$ is non-degenerate with respect to the second variable, we have that
$$^\perp T(X)=\left\{ y\in Y\mid A(x,y)=0,\text{ for all } x\in X\right\}=\left\{ 0\right\}.$$
Hence 
$$(^\perp T(X))^\perp=Y^{*}$$
and so, by \cite[Proposition 2.6.6]{Meg},
$$\overline{T(X)}^{w^{*}}=Y^{*}.$$
Thus to show that $T$ maps $X$ onto $Y^{*}$ we need to prove that $T(X)$ is $w^{*}$-closed in $Y^{*}$. To see that let
$\{ Tx_{\lambda}\}_{\lambda\in\Lambda}$ be a net in $T(X)$ and $y^{*}$ be an element of $Y^{*}$ such that
\begin{equation*}
Tx_{\lambda}\stackrel{w^*}{\longrightarrow}y^{*}.
\end{equation*}
Without loss of generality we may assume, using the special case of the Krein-\v{S}mulian Theorem on $w^{*}$-closed linear subspaces (see \cite[Corollary 2.7.12]{Meg}), the proof of which is originally due to Banach \cite[Theorem 5, p. 124]{Banach} for the separable case and to Dieudonn\'{e} \cite[Theorem 23]{Dieudonne} for the general case, that $\{ Tx_{\lambda}\}_{\lambda\in\Lambda}$ is bounded. Thus, since $\| Tx\|\geq c\| x\|$, for all $x\in X$, the net $\{ x_{\lambda}\}_{\lambda\in\Lambda}$ is also bounded. Hence, since $X$ is reflexive, there exist a subnet $\{ x_{\lambda_{\mu}}\}_{\mu\in M}$ and an element $x$ of $X$ such that $\{ x_{\lambda_{\mu}}\}_{\mu\in M}$ converges weakly to $x$. Since $T$ is $w-w^{*}$ continuous $Tx_{\lambda_{\mu}}\stackrel{w^*}{\longrightarrow}Tx$. Hence $Tx=y^{*}$ and so $T(X)$ is $w^{*}$-closed.
\end{proof}
%%%%%%%%%%%%%%%%%%%%%%%%%%%%%%%%%%%%%%%%%%%%%%%%%%%%%%%%%%%%%%%%%%%%%%%%

\begin{remark}
An alternative proof of the previous lemma can be obtained using the Closed Range Theorem.
\end{remark}
%%%%%%%%%%%%%%%%%%%%%%%%%%%%%%%%%%%%%%%%%%%%%%%%%%%%%%%%%%%%%%%%%%%%%%%%

We are now in a position to prove our main result for the linear case.
%%%%%%%%%%%%%%%%%%%%%%%%%%%%%%%%%%%%%%%%%%%%%%%%%%%%%%%%%%%%%%%%%%%%%%%%

\begin{theorem}
\label{theorem3}
Let $X$ be a reflexive Banach space, $Y$ be a Banach space, $\Lambda$ be a directed set, $\{ X_{\lambda}\}_{\lambda\in\Lambda}$ be a family of closed subspaces of $X$, $\{ Y_{\lambda}\}_{\lambda\in\Lambda}$ be an upwards directed family of closed subspaces of $Y$ and $V=\bigcup _{\lambda\in\Lambda}Y_{\lambda}$. Suppose that
\begin{equation*}
A:X\times V\longrightarrow \mathbb R
\end{equation*}
is a function for which the following hold:
\begin{itemize}
\item[\textit{(a)}] $A_{\lambda}=A\left| _{X_{\lambda}\times Y_{\lambda}}\right.$ is a bounded bilinear form, for all $\lambda\in\Lambda$.
\item[\textit{(b)}] $A(\cdot ,v)$ is a bounded linear functional on $X$, for all $v\in V$.
\item[\textit{(c)}] $A_{\lambda}$ is non-degenerate with respect to the second variable, for all $\lambda\in\Lambda$.
\item[\textit{(d)}] There exists $c>0$ such that, for all $\lambda\in\Lambda$,
\begin{equation*}
\sup_{y\in Y_{\lambda}, \| y\| =1}|A_{\lambda}(x,y)|\geq c\| x\| ,
\end{equation*}
for all $x\in X_{\lambda}$.
\end{itemize}
Then, for each bounded linear functional $v^{*}$ on $V$, there exists $x\in X$ such that
\begin{equation*}
A(x,v)=\langle v^{*},v\rangle ,
\end{equation*}
for all $v\in V$.
\end{theorem}
%%%%%%%%%%%%%%%%%%%%%%%%%%%%%%%%%%%%%%%%%%%%%%%%%%%%%%%%%%%%%%%%%%%%%%%%

\begin{proof}
Let $v^{*}\in V^{*}$ and, for each $\lambda\in\Lambda$, let $v^{*}_{\lambda}=v^{*}\left|_{Y_{\lambda}}\right.$. For all $\lambda\in\Lambda$, $v^{*}_{\lambda}$ is a bounded linear functional on $Y_{\lambda}$. By hypothesis, for all $\lambda\in\Lambda$, $A_{\lambda}$ is a bounded bilinear form on $X_{\lambda}\times Y_{\lambda}$ satisfying the two conditions of Lemma \ref{lemma1}. Since, for all $\lambda\in\Lambda$, $X_{\lambda}$ is a reflexive Banach space, we get that for each $\lambda\in\Lambda$ there exists a unique $x_{\lambda}$ such that $A_{\lambda}(x_{\lambda},y)=\langle v^{*}_{\lambda},y\rangle$, for all $y\in Y_{\lambda}$. Since $A$ satisfies condition (d), we get that, for all $\lambda\in\Lambda$,
\begin{equation*}
%\begin{array}{rcl}
c\| x_{\lambda}\|
\leq
\sup_{y\in Y_{\lambda}, \| y\| =1}|A_{\lambda}(x_{\lambda},y)|
=
\sup_{y\in Y_{\lambda}, \| y\| =1}|\langle v^{*}_{\lambda},y\rangle |
\leq
\| v^{*}\| .
%\end{array}
\end{equation*}
So $\{ x_{\lambda}\}_{\lambda\in\Lambda}$ is a bounded net in $X$. Since $X$ is reflexive, there exist a subnet $\{ x_{\lambda_{\mu}}\}_{\mu\in M}$ of $\{ x_{\lambda}\}_{\lambda\in\Lambda}$ and $x$ in $X$ such that $\{ x_{\lambda_{\mu}}\}_{\mu\in M}$ converges weakly to $x$.

We are going to prove that $A(x,v)=\langle v^{*},v\rangle$, for all $v\in V$. Take $v\in V$. Then there exists some $\lambda_{0}\in\Lambda$ with $v\in Y_{\lambda_{0}}$. Since  $\{ x_{\lambda_{\mu}}\}_{\mu\in M}$ is a subnet of $\{ x_{\lambda}\}_{\lambda\in\Lambda}$, there exists some $\mu_{0}\in M$ with $\lambda_{\mu_{0}}\geq\lambda_{0}$. 
Hence, since the family $\{ Y_{\lambda}\}_{\lambda\in\Lambda}$ is upwards directed,
\begin{equation*}
v\in Y_{\lambda_{\mu}},
\end{equation*}
for all $\mu\geq\mu_{0}$. Thus, for all $\mu\geq\mu_{0}$,
\begin{equation*}
A_{\lambda_{\mu}}(x_{\lambda_{\mu}},v)=\langle v^{*}_{\lambda_{\mu}},v\rangle.
\end{equation*}
Therefore
\begin{equation*}
\lim_{\mu\in M}A(x_{\lambda_{\mu}},v)=
\langle v^{*},v\rangle .
\end{equation*}
Since $A(\cdot ,v)$ is a bounded linear functional on $X$
\begin{equation*}
\lim_{\mu\in M}A(x_{\lambda_{\mu}},v)=A(x,v).
\end{equation*}
Hence $A(x,v)=\langle v^{*},v\rangle$.
\end{proof}
%%%%%%%%%%%%%%%%%%%%%%%%%%%%%%%%%%%%%%%%%%%%%%%%%%%%%%%%%%%%%%%%%%%%%%%%
%%%%%%%%%%%%%%%%%%%%%%%%%%%%%%%%%%%%%%%%%%%%%%%%%%%%%%%%%%%%%%%%%%%%%%%%

The following example illustrates the possible applicability of Theorem \ref{theorem3}.
%%%%%%%%%%%%%%%%%%%%%%%%%%%%%%%%%%%%%%%%%%%%%%%%%%%%%%%%%%%%%%%%%%%%%%%%

\begin{example}
\label{example1}
Let $a\in C^{1}(0,1)$ be a decreasing function with $\displaystyle\lim_{t\rightarrow 0}a(t)=\infty$ and $a(t)\geq 0$,
for all $t\in (0,1)$. We will establish the existence of a solution for the following
Cauchy problem:
\begin{equation}
\label{problem1}
\left\{
\begin{array}{rl}
u'+a(t)u=f&\text{a.e.\ on }(0,1)\\
u(0)=0&
\end{array}
\right.
\end{equation}
where $f\in L^{2}(0,1)$.

Let $X=\{ u\in H^{1}(0,1)\mid u(0)=0\}$ equipped with the norm
$\|u\| = \left(\int_0^1 |u'|^2dt\right)^\frac{1}{2}$, which is equivalent to the original
Sobolev norm, and $Y=L^{2}(0,1)$. Note that $X$ is a reflexive Banach space,
being a closed subspace of $H^{1}(0,1)$. Let $\{\alpha_{n}\}_{n\in\mathbb N}$ be a decreasing sequence in $(0,1)$
with $\displaystyle\lim_{n\rightarrow\infty}\alpha_{n}=0$. Define 
\begin{equation*}
X_{n}=\{ u\in H^{1}(\alpha_{n},1)\mid
u(\alpha_{n})=0\} \text{ and }Y_n=L^2(\alpha_{n},1)
\end{equation*}
(we can consider $X_{n}$ and $Y_n$ as  closed subspaces of $X$ and $Y$ respectively, 
by extending their elements by zero outside $(\alpha_{n},1)$). 
Also let $V=\bigcup_{n=1}^\infty Y_n$.

Let $A:X\times V\longrightarrow \mathbb R$ be the bilinear map
defined by
\begin{equation*}
A(u,v)=\int_{0}^{1}u'vdt+\int_{0}^{1}a(t)uvdt.
\end{equation*}
$A$ is well-defined and $A(\cdot,v)$ is a bounded linear functional on $X$ for any $v\in
V$.

Let $A_{n}=A|_{X_{n}\times Y_{n}}$. $A_{n}$ is a bounded, bilinear form
since
\begin{equation*}
|A_n(u,v)|\leq (1+M_n) \| u\|_{X_{n}}\| v\|_{Y_{n}}
\end{equation*}
where $M_{n}$ is the bound of $a$ on $[\alpha_{n},1]$. It should be noted that $A$ is not 
bounded on the whole of $X\times V$.

To show that $A_{n}$ is non-degenerate let $v\in Y_{n}$ and assume
that $A_{n}(u,v)=0$, for all $u\in X_{n}$, i.e.
\begin{equation*}
\int_{\alpha_{n}}^{1}(u'+a(t)u)vdt=0,\text{ for all }u\in X_{n}.
\end{equation*}
It is easy to see that the above implies that
\begin{equation*}
\int_{\alpha_{n}}^{1}wvdt=0 ,
\end{equation*}
for any continuous function $w$ and therefore $v=0$.

We next show that
\begin{equation*}
\sup_{\| v\| =1,\;v\in Y_{n}}|A_{n}(u,v)|\geq \|u\|_{X_{n}}.
\end{equation*}
Define $T_{n}:X_{n}\longrightarrow Y_{n}^*$ by $\langle T_{n}u,v\rangle =A_{n}(u,v)$.
$T_n$ is a well-defined bounded linear operator and $T_{n}u=u'+a(t)u$. Hence
\begin{equation*}
\begin{array}{rcl}
\| T_{n}u\|^{2}
&=&
\displaystyle{\int_{\alpha_{n}}^{1}|u'+a(t)u|^2dt}\\[15pt]
&=&\displaystyle{\int_{\alpha_{n}}^{1}|u'|^2dt+\int_{\alpha_{n}}^{1}a^{2}(t)|u|^2dt+\int_{\alpha_{n}}^{1}a(t)(u^2)'dt}\\[15pt]
&=&\displaystyle{\int_{\alpha_{n}}^{1}|u'|^2dt+\int_{\alpha_{n}}^{1}(a^{2}(t)-a'(t))|u|^2dt+a(1)u^{2}(1)}\\[15pt]
&\geq&
\| u\|^{2}_{X_{n}},
\end{array}
\end{equation*}
since $u(\alpha_n)=0$, $a$ is decreasing and $a(t)\geq 0$, for all $t\in (0,1)$.

All the hypotheses of Theorem \ref{theorem3} are hence satisfied and so if $F\in
V^{*}$ is defined by $F(v)=\int_{0}^{1}fvdt$, then there exists $u\in X$ such that
\begin{equation*}
A(u,v)=F(v) ,\text{ for all }v\in V.
\end{equation*}
Thus $u$ satisfies (\ref{problem1}).
\end{example}
%%%%%%%%%%%%%%%%%%%%%%%%%%%%%%%%%%%%%%%%%%%%%%%%%%%%%%%%%%%%%%%%%%%%%%%%%%%%%%%%%%

\section{The nonlinear case}
%%%%%%%%%%%%%%%%%%%%%%%%%%%%%%%%%%%%%%%%%%%%%%%%%%%%%%%%%%%%%%%%%%%%%%%%

We start by recalling some well-known definitions:
%%%%%%%%%%%%%%%%%%%%%%%%%%%%%%%%%%%%%%%%%%%%%%%%%%%%%%%%%%%%%%%%%%%%%%%%

\begin{definition}
\label{definition1}
Let $T:X\longrightarrow X^{*}$ be an operator. We say that $T$ is:
\begin{itemize}
\item[(i)] monotone if $\langle Tx-Ty,x-y\rangle\geq 0$, for all
$x,\;y\in X$.
\item[(ii)] hemicontinuous if, for all $x,\;y\in X$, $T(x+ty)\stackrel{w}{\longrightarrow}Tx$ as
$t\longrightarrow 0^{+}$.
\item[(iii)] coercive if
\begin{equation*}
\lim_{||x||\rightarrow \infty}\frac{\langle Tx,x\rangle}{||x||}=\infty.
\end{equation*}
\end{itemize}
\end{definition}
%%%%%%%%%%%%%%%%%%%%%%%%%%%%%%%%%%%%%%%%%%%%%%%%%%%%%%%%%%%%%%%%%%%%%%%%

We also need the following generalization of the notion of type $M$ operator (for the classical definition see \cite {Brezis} or \cite{Zeidler}).
%%%%%%%%%%%%%%%%%%%%%%%%%%%%%%%%%%%%%%%%%%%%%%%%%%%%%%%%%%%%%%%%%%%%%%%%

\begin{definition}
Let $X$ be a Banach space, $V$ be a linear subspace of $X$ and 
\begin{equation*}
A:X\times V\longrightarrow \mathbb{R}
\end{equation*}
be a function. We say that $A$ is of type $M$
with respect to $V$ if, for any net $\{ v_\lambda\}_{\lambda\in\Lambda}$ in $V$, $x\in X$ and $v^{*}\in V^*$,
\begin{itemize}
\item[\textit{(a)}] $v_\lambda\stackrel{w}{\longrightarrow}x$
\item[\textit{(b)}] $A(v_\lambda,v)\longrightarrow\langle v^{*},v\rangle$, for all $v\in V$
\item[\textit{(c)}] $A(v_\lambda,v_\lambda)\longrightarrow\langle \hat{v}^{*},x\rangle$, where $\hat{v}^{*}$ is the extension of $v^{*}$ on the closure of $V$
\end{itemize}
imply that $A(x,v)=\langle v^{*},v\rangle$, for all $v\in V$.
\end{definition}
%%%%%%%%%%%%%%%%%%%%%%%%%%%%%%%%%%%%%%%%%%%%%%%%%%%%%%%%%%%%%%%%%%%%%%%%

Our result is the following:
%%%%%%%%%%%%%%%%%%%%%%%%%%%%%%%%%%%%%%%%%%%%%%%%%%%%%%%%%%%%%%%%%%%%%%%%

\begin{theorem}
\label{theorem4}
Let $X$ be a reflexive Banach space, $\Lambda$ be a directed set,
$\{ X_{\lambda}\}_{\lambda\in\Lambda}$ be an upwards directed family of closed subspaces of $X$ and
$V=\bigcup _{\lambda\in\Lambda}X_{\lambda}$. Suppose that 
\begin{equation*}
A:X\times V\longrightarrow \mathbb{R}
\end{equation*} 
is a function for which the following hold:
\begin{itemize}
\item[\textit{(a)}] $A$ is of type M with respect to $V$.
\item[\textit{(b)}] $\displaystyle{\lim_{||x||\rightarrow \infty}\frac{A(x,x)}{||x||}=\infty}$.
\item[\textit{(c)}] $A_{\lambda}(x,\cdot)\in X_{\lambda}^{*}$,  for all $\lambda\in\Lambda$ and all $x\in X_{\lambda}$, where $A_{\lambda}$ is the restriction of $A$ on $X_{\lambda}\times X_{\lambda}$.
\item[\textit{(d)}] The operator $T_{\lambda}:X_{\lambda}\longrightarrow X_{\lambda}^*$, defined by
$\langle T_{\lambda}x,y\rangle=A_{\lambda}(x,y)$, for all $x,y\in X_{\lambda}$, is
monotone and hemicontinuous, for all $\lambda\in\Lambda$.
\end{itemize}
Then for each $v^{*}\in V^{*}$ there exists $x\in X$ such that
\begin{equation*}
A(x,v)=\langle v^{*},v\rangle,
\end{equation*}
for all $v\in V$.
\end{theorem}
%%%%%%%%%%%%%%%%%%%%%%%%%%%%%%%%%%%%%%%%%%%%%%%%%%%%%%%%%%%%%%%%%%%%%%%%

\begin{proof}
As in the proof of Theorem \ref{theorem3}, for each $\lambda\in\Lambda$, let $v^{*}_{\lambda}=v^{*}\left|_{X_{\lambda}}\right.$.
By the Browder-Minty Theorem (see \cite[Theorem 26.A]{Zeidler}), a monotone, coercive and hemicontinuous operator, from a real reflexive Banach space into its dual, is onto.
Thus, by (b) and (d), for each $\lambda\in\Lambda$, the operator $T_{\lambda}$ is onto and so
there exists $x_{\lambda}\in X_{\lambda}$ such that
\begin{equation*}
A_{\lambda}(x_{\lambda},y)=\langle v^{*}_{\lambda},y\rangle,
\end{equation*}
for all $y\in X_{\lambda}$. In particular
$A_{\lambda}(x_{\lambda},x_{\lambda})=\langle v^{*}_{\lambda},x_{\lambda}\rangle$ and hence, by (b), we get that the net $\{ x_{\lambda}\}_{\lambda\in\Lambda}$ is bounded. Continuing as in the proof of Theorem \ref{theorem3} and applying the fact that $A$ is of type $M$ with respect to $V$ we get the required result.
\end{proof}
%%%%%%%%%%%%%%%%%%%%%%%%%%%%%%%%%%%%%%%%%%%%%%%%%%%%%%%%%%%%%%%%%%%%%%%%

\begin{remark}
It should be noted that, since a crucial point in the above proof is the existence and boundedness of the net $\{ x_{\lambda}\}_{\lambda\in\Lambda}$, variants of the previous theorem could be obtained using in (b) and (d) alternative conditions corresponding to other surjectivity results.
\end{remark}
%%%%%%%%%%%%%%%%%%%%%%%%%%%%%%%%%%%%%%%%%%%%%%%%%%%%%%%%%%%%%%%%%%%%%%%%

We now apply Theorem \ref{theorem4} to a singular Dirichlet problem.
%%%%%%%%%%%%%%%%%%%%%%%%%%%%%%%%%%%%%%%%%%%%%%%%%%%%%%%%%%%%%%%%%%%%%%%%
\begin{example}
\label{example2}
Let $\Omega$ be a bounded domain in $\mathbb R^N$. We consider the Dirichlet problem
\begin{equation}
\label{equation2}
\left\{
\begin{array}{rl}
\displaystyle-\sum_{i=1}^{N}\frac{\partial}{\partial x_i}\left(a(x)\frac{\partial u}{\partial x_i}\right)+f(x,u)=0 & \text{a.e.\ on }\Omega\\[5pt]
u=0&\text{ on } \partial \Omega
\end{array}
\right.,
\end{equation}
where $a\in L^{\infty}_{loc}(\Omega)$ and there exists $c_1>0$ such that $a(x)\geq c_1$ a.e.\ on $\Omega$ and $f:\Omega\times\mathbb R\longrightarrow \mathbb R$ is a monotone increasing (with respect to its second variable for each fixed $x\in \Omega$), Carath\'{e}odory function, for which there exist $h\in L^2(\Omega)$ and $c_2>0$ such that
%%%%%%%%%%%%%%%%%%%%%%%%%%%%%%%%%%%%%%%%%
\begin{equation}
\label{growth}
|f(x,u)|\leq h(x)+c_2|u|, \text{ for all } x\in\Omega \text{ and } u\in\mathbb R.
\end{equation}
%%%%%%%%%%%%%%%%%%%%%%%%%%%%%
We will show that if the above hypotheses on $a$ and $f$ hold, then problem (\ref{equation2}) has a weak solution, i.e. that there exists a function $u\in H^1_0(\Omega)$ with
$$\int_\Omega a(x)\nabla u\nabla v dx+\int_\Omega f(x,u)v dx=0, \text{ for all }v\in C^\infty_0(\Omega).$$

To this end let $X=H^{1}_{0}(\Omega)$, $\{\Omega_n\}_{n\in\mathbb{N}}$ be an increasing sequence of open subsets of $\Omega$ such that $\overline{\Omega_n}\subseteq\Omega_{n+1}$ and 
$$\bigcup_{n=1}^{\infty}\Omega_n=\Omega$$ 
and $X_{n}=H^{1}_{0}(\Omega_n)$, for each $n\in\mathbb{N}$. Observe that we can consider each $X_n$ as a closed subspace of $X$ by extending its elements by zero outside $\Omega_n$ and let 
$$V=\bigcup_{n=1}^\infty X_n.$$ 
Finally let
\begin{equation*}
A:X\times V\longrightarrow \mathbb{R}
\end{equation*} 
be the function defined by
\begin{equation*}
A(u,v)=\int_\Omega a(x)\nabla u\nabla v dx+\int_\Omega f(x,u)v dx.
\end{equation*}
%%%%%%%%%%%%%%%%%%%%%%%%%%%%%%%

By $a(x)\geq c_1$ a.e.\ on $\Omega$, the monotonicity of $f$ and the growth condition (\ref{growth}), we have
%%%%%%%%%%%%%%%%%%%%%%%%%
\begin{eqnarray*}
\label{coercive}
A(u,u)&=&\int_\Omega a(x)|\nabla u|^2 dx+\int_\Omega f(x,u)u dx\\[10pt]
&=&\int_\Omega a(x)|\nabla u|^2 dx+\int_\Omega (f(x,u)-f(x,0))u dx+\int_\Omega f(x,0)u dx\\[10pt]
&\geq & c_1 \|\nabla u\|^2_{L^2(\Omega)}-\| h\|_{L^2(\Omega)}\| u\|_{H^1_0(\Omega)}.
\end{eqnarray*}
%%%%%%%%%%%%%%%%%%%%%%%%%
Since by the Poincar\'{e} inequality $\| \nabla u\|_{L^2(\Omega)}$ is equivalent to the norm of $X$ it follows that $A$ is coercive.

Let $A_n=A|_{X_n\times X_n}$. Then, since $a\in L^\infty_{loc}(\Omega)$ it follows that $a\in L^\infty(\Omega_n)$, for all $n\in \mathbb N$. Combining this with (\ref{growth}), we have that
$$|A_n(u,v)|\leq c(u,n) \|v\|_{X_n},$$
where $c(u,n)$ is a positive constant depending on $n$ and $u$. So the operator 
$$T_n:X_n\longrightarrow X_n^\ast,$$ with $\langle T_nu,v\rangle_{X_n}=A_n(u,v)$ is well-defined, for all $n\in \mathbb N$. Let 
$$T_{1,n},\,T_{2,n}:X_n\longrightarrow X_n^\ast$$ 
be the operators defined by
$$\langle T_{1,n}u,v\rangle_{X_n}=\int_{\Omega_n}a(x)\nabla u\nabla v dx\;\text{ and }\;\langle T_{2,n}u,v\rangle_{X_n}=\int_{\Omega_n} f(x,u) v dx.$$
Then $T_{1,n}$ is a monotone bounded linear operator. Using the monotonicity of $f$, it is easy to see that $T_{2,n}$ is monotone. Finally, recalling that the Nemytskii operator corresponding to $f$ is continuous (see for example \cite[Proposition 26.7]{Zeidler}) and that the embedding of $X_n$ into $L^2(\Omega_n)$ is compact we have that $T_{2,n}$ is hemicontinuous. Thus $T_n=T_{1,n}+T_{2,n}$ is monotone and hemicontinuous, for all $n\in \mathbb N$. 

To finish the proof let $u_n\stackrel{w}{\longrightarrow}u$ in $X$. Then since, for all $v\in V$,
$$
u\longmapsto\int_\Omega a(x)\nabla u\nabla v dx
$$
is a bounded linear functional and, by the continuity of the Nemytskii operator and the compactness of the embedding of $X$ into $L^2(\Omega)$, 
$$\int_\Omega f(x,u_n)v dx\longrightarrow\int_\Omega f(x,u)v dx,$$
for all $v\in V$, we get that
$$A(u_n,v)\longrightarrow A(u,v), \text{ for all } v\in V.$$ 
Thus $A$ is of type $M$ with respect to $V$. Applying now Theorem \ref{theorem4} we get that there exists $u\in X$ such that $A(u,v)=0$, for all $v\in V$. Observing that $C^\infty_0(\Omega)$ is contained in $V$ we get that $u$ is the required weak solution of (\ref{equation2}).
\end{example}
%%%%%%%%%%%%%%%%%%%%%%%%%%%%%%%%%%%%%%%%%%%%%%%%%%%%%%%%%%%%%%%%%%%%%%%%%

\smallskip
%%%%%%%%%%%%%%%%%%%%%%%%%%%%%%%%%%%%%%%%%%%%%%%%%%%%%%%%%%%%%%%%%%%%%%%%
\noindent {\it Acknowledgments.} The authors would like to thank Prof. A. Katavolos for pointing out an error in an earlier version of this paper and the two referees for their comments and suggestions which improved both the content and the presentation of this paper.
%%%%%%%%%%%%%%%%%%%%%%%%%%%%%%%%%%%%%%%%%%%%%%%%%%%%%%%%%%%%%%%%%%%%%%%%

\end{document}